\documentclass[12pt]{sagej}
\usepackage{amsmath,amsfonts,amssymb,amsthm,epsfig,subfigure,graphicx,url}
\usepackage{cancel}
\usepackage{eufrak}
\usepackage[utf8]{inputenc}

\newcommand{\be}{\begin{equation}}
\newcommand{\en}{\end{equation}}

\renewcommand{\vec}[1]{\boldsymbol{#1}}

\begin{document}

\runninghead{Destrade and Saccomandi}



\title{Plane-polarised finite-amplitude shear waves in deformed incompressible materials}


\author{ Michel Destrade\affilnum{1,2} and  Giuseppe Saccomandi\affilnum{3,1}}
\affiliation
{\affilnum{1}School of Mathematical and Statistical Sciences,  
NUI Galway, Ireland. \\
\affilnum{2}Department of Engineering Mechanics, Zhejiang University, Hangzhou 310027, PR China.
\\
\affilnum{3}Dipartimento di Ingegneria, 
  Universit\`{a} degli Studi di Perugia and Sezione INFN of Perugia, Italy.}

\corrauth{Michel Destrade, School of Mathematical and Statistical Sciences,  
NUI Galway. }
\email{michel.destrade@nuigalway.ie}


\begin{abstract}

We investigate how two finite-amplitude, transverse, plane body waves may be superposed to propagate in a deformed hyperelastic incompressible solid.
We find that the equations of motion reduce to a well-determined system of partial differential equations, making the motion controllable for all solids. 
We find that in deformed Mooney-Rivlin materials, they may travel along any direction and be polarised along any transverse direction, an extension of a result by Boulanger and Hayes [Quart. J. Mech. Appl. Math. 45 (1992) 575]. 
Furthermore, their motion is governed by a linear system of partial differential equations, making the Mooney-Rivlin special in that respect. 
We select another model to show that for other materials, the equations are nonlinear.
We use asymptotic equations to reveal the onset of nonlinearity for the waves, paying particular attention to how close the propagation direction is to the principal axes of pre-deformation. 
\\[24pt]
\emph{Dedicated to the late Hui-Hui Dai, a true gentleman, always kind, sincere and generous; greatly missed.}


\end{abstract}

\maketitle



\section{Introduction}


In 1969, Currie and Hayes \cite{Currie} showed that two linearly polarised finite-amplitude transverse waves, polarised in two orthogonal directions, may propagate along any direction in a Mooney-Rivlin material maintained in a state of arbitrary static finite homogeneous deformation. 
Their results were a generalisation, to some extent, of previous findings contained by Green \cite{Green} and Carroll \cite{Caroll}.
Later,  beginning with a  paper published  in 1992 \cite{BH1992}, Boulanger and Hayes wrote a series of papers to investigate and extend those results further, in a deep and elegant manner; their findings were collected and summarised in the survey \cite{BH2001}.
 
Boulanger and Hayes's main finding is that if $\vec{B}$ denotes the left Cauchy-Green deformation tensor of the homogeneous pre-deformation and $\vec{n}$ denotes the direction of propagation, then the directions of polarisation of the two possible transverse waves propagating along $\vec n$ must be along the principal axes of \color{black} the elliptical section of the  $\vec{x \cdot B}^{-1} \vec{x}=1$ ellipsoid by the $\vec{x\cdot n = 0}$ plane. \color{black}
That restriction is in place for the class of Mooney-Rivlin materials, that is, those incompressible, homogeneous, isotropic materials with strain energy density of the form $W= C(I_1-3) + E(I_2-3)$, where $I_1$ and $I_2$ are the first and second principal invariants of $\vec C$, the right Cauchy-Green deformation tensor, and $C$, $E$ are positive constants.  
Only in the special case of a neo-Hookean material, when $C_2=0$, may the directions of polarisation be along any direction orthogonal to $\vec{n}$.

This result is the natural consequence  of the \emph{propagation condition},
\be \label{i1}
\vec{a  \cdot B}^{-1} \vec{b}=0,
\en 
where ($\vec{n}, \vec{a}, \vec{b}$) is an orthonormal triad, with $\vec{n}$ the direction of propagation.
This condition is satisfied when the vectors $\vec a$ and $\vec b$ are aligned with the axes of the central elliptic section of the ellipsoid by the plane $\vec{n\cdot x}=0$. 
Then two single transverse waves may propagate individually, one polarised along $\vec{a}$, the other along $\vec{b}$. 
The propagation condition \eqref{i1} follows from the equations of motion: these turn out to be an overdetermined system of three partial differential equations for two unknowns functions, the pressure field and the transverse linearly-polarized wave function. To reduce the system to two equations for two unknowns, the propagation condition \eqref{i1} must apply.

Another result established by Boulanger and Hayes \cite{BH1992} concerns the \emph{superposition} of two waves which propagate along the same direction $\vec{n}$. Indeed, they found that for any propagation direction $\vec{n}$, \emph{two} waves could propagate and solve the equations of motion simultaneously if they were linearly polarised along the unit vectors $\vec a$ and $\vec b$ solutions to Equation \eqref{i1}. This is a remarkable result, because no assumption is made about the magnitude of the waves, and the theory is completely non-linear and exact.

In the present note we show that if we consider two shear waves propagating along \emph{any} direction $\vec{n}$  and polarized along \emph{any} orthogonal unit vectors $\vec{a}$ and $\vec{b}$ in the $\vec{n \cdot x} = 0$ plane (not necessarily satisfying \eqref{i1}), in \emph{any} deformed  hyperelastic material (not just the Mooney-Rivlin materials), then the equations of motion reduce to a well-determined system of partial differential equations. 
It follows that this motion is controllable for any incompressible and isotropic hyperelastic material.  
This is shown in the next section.

Further, we find in Section \ref{Section3} that in the special case of a Mooney-Rivlin material, the determining equations are linear and they possess solutions of permanent form, similar to what was seen in \cite{BH1992}, but here for shear waves polarised along \emph{any} orthogonal unit vectors $\vec a$ and $\vec b$.
 
Finally in Section \ref{Section4}, we take a specific form of the strain energy density to study an example of the nonlinear equations generated when the solid is not as special as the one modelled by the Mooney-Rivlin material.
We use asymptotic expansions in the amplitude to reveal the onset of nonlinearity in the equations of motion, and how great care must be taken when the direction of propagation is close to, or along a principal axis of pre-deformation. 
    

\section{Superposition of two shear  waves}

 
Let $\vec{X}$ be the position vector of a particle in a hyperelastic body in the reference configuration, and $\vec{x}$ be its position vector in the current configuration. 
We call $\vec{F}=\partial \vec{x}/\partial \vec{X}$  the deformation gradient and $\vec{B}=\vec{F}\vec{F}^T$, $\vec C = \vec F^T \vec F$ the left and right Cauchy-Green deformation tensors, respectively.

We consider homogeneous, incompressible ($\det \vec{F}=1$ at all times), iso\-tro\-pic and hyperelastic materials, so that their strain-energy density $W$ (measured per unit volume in the undeformed state) is a function of the form
\begin{equation}
W=W(I_1,I_2),
\end{equation}
where $I_1=\text{tr}\,\vec{C}$ and $I_2=\text{tr}\, \vec{C}^{-1}$ are the first and second principal invariants of $\vec{C}$, respectively. 
Then the Cauchy stress $\vec{T}$ is 
\be \label{e1}
\vec{T}=-p\vec{I}+2W_1\vec{B}-2W_2\vec{B}^{-1},
\en
where $p$ is the Lagrange multiplier associated with the constraint of incompressibility, and $W_k = \partial W/\partial I_k$.
The equations of motion, in the absence of body forces, are
\be \label{e2}
\text{div} \, \vec{T} = \rho \; \partial^2 \vec x/\partial t^2,
\en
 where $\rho$ is the mass density (which is constant, because of incompressibility).
 
First we consider that the solid is subject to a static finite homogeneous deformation:
\be
\overline{\vec{x}}=\overline{\vec{F}} \vec{X},
\en
where $\overline{\vec{F}}$ is a constant tensor, such that $\det \overline{\vec{F}}=1$ (to accommodate the constraint of incompressibility).  
Because this deformation is universal, the  stress tensor required to support it is, according to \eqref{e1},
\be \label{e3}
\overline{\vec{T}}=-\overline{p}\vec{I}+2W_1(\overline{I}_1, \overline{I}_2)\overline{\vec{B}}-2W_2(\overline{I}_1, \overline{I}_2)\overline{\vec{B}}^{-1}.
\en
Taking the Lagrange multiplier $\overline p$ to be a constant ensures that the equations of static equilibrium $\text{div} \, \overline{\vec T} = \vec 0$ are satisfied.

Then we superimpose on this deformation two plane homogeneous body shear waves, both propagating along the direction of the unit vector $\vec{n}$, with one polarised along $\vec{a}$ (a unit vector orthogonal to $\vec n$) and one along $\vec{b}=\vec{n} \times \vec{a}$.
Hence this motion is of the form
\be \label{e4}
\vec{x}=\overline{\vec{x}}+f(\overline{\eta}, t) \vec{a}+g(\overline{\eta}, t) \vec{b}, \qquad p=\overline{p}+q(\overline{\eta}, t),
\en
where $(\overline{\eta}, \overline{\xi}, \overline{\zeta})$ are the components of $\overline{\vec{x}}$ in the ($\vec{n}, \vec{a}, \vec{b}$) orthonormal basis, and $f,g,q$ are yet unknown amplitude functions. 

The deformation gradient associated with this motion is 
\be
\vec{F}=\left(\vec{I}+f_{\overline{\eta}} \vec{a} \otimes \vec{n}+g_{\overline{\eta}} \vec{b} \otimes \vec{n}\right)\overline{\vec{F}},
\en
where the subscript denotes partial differentiation.
Clearly, this is an isochoric motion, respecting the constraint of incompressibility.

Let $(\eta, \, \xi , \, \zeta)$ be the components of the position vector $\vec{x}$ in ($\vec{n}, \vec{a}, \vec{b}$).
Then the motion \eqref{e4}$_1$ reads
\be \label{e5}
\eta=\overline{\eta}, \quad \xi=\overline{\xi}+f(\eta,t), \quad \zeta=\overline{\zeta}+g(\eta,t).
\en 
We then compute the following kinematic quantities associated with the motion,
\begin{align} \notag \label{e6}
& \vec{B}=\left( \vec{I}+f_\eta \vec{a}\otimes\vec{n}+g_\eta \vec{b}\otimes \vec{n} \right)\overline{\vec{B}}\left( \vec{I}+f_\eta \vec{n}\otimes\vec{a}+g_\eta \vec{n}\otimes \vec{b} \right),
 \notag \\
& \vec{B}^{-1}=\left( \vec{I}-f_\eta \vec{n}\otimes\vec{a}-g_\eta \vec{n}\otimes \vec{b} \right)\overline{\vec{B}}^{-1}\left( \vec{I}-f_\eta \vec{a}\otimes\vec{n}-g_\eta \vec{b}\otimes \vec{n} \right),
 \notag \\
& I_1=\overline{I}_1+2(f_\eta \vec{n \, \cdot} \overline{\vec{B}} \vec{a}+g_\eta \vec{n \, \cdot} \overline{\vec{B}} \vec{b})+\left(f^2_\eta+g^2_\eta \right) \vec{n \, \cdot} \overline{\vec{B}} \vec{n}, \notag \\
&  I_2=\overline{I}_2-2\left(f_\eta \vec{n \, \cdot} \overline{\vec{B}}^{-1} \vec{a}+g_\eta \vec{n \, \cdot} \overline{\vec{B}}^{-1} \vec{b}\right)+ f^2_\eta \vec{a \,\cdot} \overline{\vec{B}}^{-1} \vec{a}+g^2_\eta\vec{b \,  \cdot} \overline{\vec{B}}^{-1} \vec{b}.
\end{align}

Now the equations of motion \eqref{e2}, written in the basis $(\vec{n}, \vec{a}, \vec{b})$, reduce to
\be \label{e7}
0=\frac{\partial T_{\eta \eta } }{\partial \eta}, \qquad \rho f_{tt}=\frac{\partial T_{\xi \eta} }{\partial  \eta}, \qquad \rho g_{tt}=\frac{\partial T_{\zeta \eta} }{\partial \eta}.
\en 
The first of these is satisfied by an appropriate choice of the Lagrange multiplier $q$, which then can be forgotten about, as it does not play any role in the other two equations.

Introducing the notation  $\mathcal{W}=W(I_1, I_2)$ where $I_1$ and $I_2$ are now given by \eqref{e6}, we compute the remaining non-zero components  of the stress tensor as
\begin{multline}  \label{t1}
 T_{\xi \eta}=2\mathcal{W}_1\left(\vec{n \, \cdot} \overline{\vec{B}}\vec{a}+f_\eta  \vec{n \, \cdot} \overline{\vec{B}}\vec{n}\right)\\
 -2\mathcal{W}_2\left(\vec{n \, \cdot} \overline{\vec{B}}^{-1}\vec{a}-f_\eta  \vec{a \,\cdot} \overline{\vec{B}}^{-1}\vec{a}-g_\eta \vec{a \, \cdot} \overline{\vec{B}}^{-1}\vec{b}\right),
\end{multline}
and
\begin{multline}  \label{t2}
T_{\zeta \eta}=2\mathcal{W}_1\left(\vec{n \, \cdot} \overline{\vec{B}}\vec{b}+g_\eta  \vec{n \, \cdot} \overline{\vec{B}}\vec{n}\right)
\\  -2\mathcal{W}_2\left(\vec{n \,\cdot} \overline{\vec{B}}^{-1}\vec{b}-f_\eta  \vec{b \, \cdot} \overline{\vec{B}}^{-1}\vec{a}-g_\eta \vec{b \,\cdot} \overline{\vec{B}}^{-1}\vec{b}\right),
\end{multline}
\color{black}
where $\mathcal W_k=\partial \mathcal W/\partial I_k$.
\color{black}

%

It is now clear that the remaining equations of motion \eqref{e7}$_{2,3}$ are a nonlinear system of two coupled differential equations for the two unknown functions $f$ and $g$ and thus, that the motion is controllable for any incompressible, isotropic, hyperelastic material. 


\section{The special case of Mooney-Rivlin materials}
\label{Section3}


Consider now the special case of  Mooney-Rivlin materials, with strain energy density 
\be \label{mr1}
W=C(I_1-3)+E(I_2-3),
\en
where $C,E$ are  material constants such that $C>0$ and $E \geq 0$. 

For these materials the system of equations \eqref{e7}$_{2,3}$, is \emph{linear}, as the equations read
\begin{align} \label{mr2}
& \rho f_{tt}= \left(C \vec{n \,\cdot} \overline{\vec{B}}\vec{n}+E \vec{a \, \cdot} \overline{\vec{B}}^{-1}\vec{a}\right)f_{\eta \eta}+ E\left(\vec{a \,\cdot} \overline{\vec{B}}^{-1}\vec{b} \right) g_{\eta \eta},  \notag
\\
& \rho g_{tt}= \left(C \vec{n \, \cdot} \overline{\vec{B}}\vec{n}+E \vec{b \,\cdot} \overline{\vec{B}}^{-1}\vec{b}\right)g_{\eta \eta}+ E \left(\vec{a \,\cdot} \overline{\vec{B}}^{-1}\vec{b} \right) f_{\eta \eta}. 
\end{align}
 
As  noted by Boulanger and Hayes \cite{BH1992}, these equations decouple  when $\vec{a}\cdot \overline{\vec{B}}^{-1}\vec{b}=0$. 
But because \eqref{mr2} is a linear system, we may in fact solve it in any case, by writing it as
\begin{equation}
\rho \begin{bmatrix} f \\ g\end{bmatrix}_{tt} = 
  \begin{bmatrix} 
 C \vec{n \,\cdot} \overline{\vec{B}}\vec{n}+E \vec{a \, \cdot} \overline{\vec{B}}^{-1}\vec{a} &   E\vec{a \,\cdot} \overline{\vec{B}}^{-1}\vec{b} \\ 
 E\vec{a \,\cdot} \overline{\vec{B}}^{-1}\vec{b} & C \vec{n \, \cdot} \overline{\vec{B}}\vec{n}+E \vec{b \,\cdot} \overline{\vec{B}}^{-1}\vec{b}
 \end{bmatrix}
 \begin{bmatrix} f \\ g\end{bmatrix}_{\eta\eta},
 \end{equation}
 and diagonalising the matrix. 
 We then arrive at the decoupled equations
 \begin{equation}
 \rho u_{tt} = \lambda_1 u_{\eta\eta}, \qquad 
 \rho v_{tt} = \lambda_2 v_{\eta\eta}, 
 \label{wave-equation}
 \end{equation}
 where the (real) eigenvalues are
 \begin{multline}
 \lambda_{1,2} = \tfrac{1}{2} 	\left[2C \vec{n \, \cdot} \overline{\vec{B}}\vec{n} + E \left(\vec{a \, \cdot} \overline{\vec{B}}^{-1}\vec{a}+ \vec{b \, \cdot} \overline{\vec{B}}^{-1}\vec{b}\right) \right. \\ 
 \left.  \pm \, E\sqrt{\left(\vec{a \,\cdot} \overline{\vec{B}}^{-1}\vec{a}-\vec{b \,\cdot} \overline{\vec{B}}^{-1}\vec{b}\right)^2+4\left(\vec{a\, \cdot} \overline{\vec{B}}^{-1}\vec{b} \right)^2}
 \right].
 \end{multline}
  and the functions $u$, $v$ are defined by
  \begin{equation}
  \begin{bmatrix}
  u \\ v \end{bmatrix}
  = \begin{bmatrix}
   \cos \theta & \sin \theta \\ -\sin\theta & \cos \theta
   \end{bmatrix}
  \begin{bmatrix}
  f \\ g \end{bmatrix},
 \quad \text{with } \tan 2\theta =  \dfrac{2 \vec{a \,\cdot} \overline{\vec{B}}^{-1}\vec{b}}
 {\vec{a \,\cdot} \overline{\vec{B}}^{-1}\vec{a}-\vec{b \,\cdot} \overline{\vec{B}}^{-1}\vec{b}}.
   \end{equation}
 Note that both eigenvalues are positive, because $\lambda_1$ is the sum of positive quantities, and writing that $\lambda_2>0$ is equivalent to 
 \begin{multline}
 C^2 \left( \vec{n \, \cdot} \overline{\vec{B}}\vec{n} \right)^2 + CE  \left( \vec{n \, \cdot} \overline{\vec{B}}\vec{n} \right)\left(\vec{a \, \cdot} \overline{\vec{B}}^{-1}\vec{a}+ \vec{b \, \cdot} \overline{\vec{B}}^{-1}\vec{b}\right)
 \\
  + E^2 \left[\left(\vec{a \,\cdot} \overline{\vec{B}}^{-1}\vec{a}\right) \left(\vec{b \,\cdot} \overline{\vec{B}}^{-1}\vec{b}\right) - \left(\vec{a\, \cdot} \overline{\vec{B}}^{-1}\vec{b} \right)^2\right]>0,
 \end{multline}
 which is always true, because the last bracketed term is $|\vec V^{-1}\vec a \times \vec V^{-1}\vec b|^2$, where $\vec V$ is the square root of $\vec B$.
 
 Both $u$ and $v$ satisfy the wave equation  \eqref{wave-equation}, with general solution
 \begin{equation}
 u = u_-(\eta - c_1t) +  u_+(\eta + c_1t), \qquad
 v = v_-(\eta - c_2t) +  v_+(\eta + c_2t), 
 \end{equation}
 where $u_\pm$, $v_\pm$ are arbitrary functions, and  the speeds are  $c_i = \sqrt{\lambda_i/\rho}$.
 
 Hence we have established that we may use Boulanger and Hayes's solution (where $\vec a$ and $\vec b$ satisfy \eqref{i1}) to generate a solution with any direction of polarization. 
 Morevover, the motion thus generated has components which travel at different speeds $c_1$ and $c_2$.

To conclude this section, we note that for Mooney-Rivlin materials, the determining equations for transverse waves are a completely exceptional hyperbolic system, once the strong ellipticity condition is ensured, see \cite{Vitolo}. 
The possibility of wave motion as solution to linear partial differential equations for the Mooney-Rivlin and neo-Hookean materials is not restricted to the case of transverse waves, as has been noticed by several authors, including Lei and Hung \cite{Lei}, Rajagopal \cite{Raja}, and Hill and Dai \cite{Hill}.


\section{Nonlinear waves}
\label{Section4}


Finite amplitude plane waves propagating in a deformed Mooney-Rivlin material lead to an exceptional hyperbolic system, as the motion is determined by a linear system of partial differential equations. 
Here we turn our attention to materials which are not special in that respect.
It leads us to a nonlinear  hyperbolic system, which we approach using asymptotic expansions. 

First we consider the \emph{case of no pre-strain}, so that $\overline{\vec{B}}\equiv\vec{I}$ and
\begin{equation}
T_{\xi \eta}=2(\mathcal{W}_1+ \mathcal{W}_2) f_\eta, \qquad T_{\zeta \eta}=2(\mathcal{W}_1 +  \mathcal{W}_2) g_\eta.
\end{equation}
Introducing the unknown functions $\hat F=f_\eta$, $\hat G=g_\eta$, we have $I_1=I_2=F^2+G^2+3$, and the remaining two determining equations in  \eqref{e7} may be rewritten as 
\be \label{nw1}
\rho \hat F_{tt}=(\mathcal{Q}\hat F)_{\eta \eta}, \qquad \rho \hat G_{tt}=(\mathcal{Q}\hat G)_{\eta \eta}, 
\en
where $\mathcal{Q} = 2 (\mathcal{W}_1 + \mathcal{W}_2)$ is the generalized shear modulus.
Clearly  $\mathcal{Q}=\mathcal{Q}(\hat F^2 + \hat G^2)$.  
Then considering that the amplitudes are \color{black} small\color{black}, we write $\hat F= \epsilon F$ and $\hat G=\epsilon G$ where $|\epsilon | \ll 1$ and expand $\mathcal{Q}$ as $\mathcal{Q}=\mu_0+\mu_1 \epsilon^2 (F^2+G^2)+\ldots$.
For the Mooney-Rivlin material, $\mu_0=2(C+E)$ and $\mu_1=\mu_2=\ldots = 0$, but in general, $\mu_1 \ne 0$ (although there are materials other than the Mooney-Rivlin material such that $\mu_1=\mu_2=\ldots = 0$, see \cite{Mangan}.)
In that case we introduce the scaled time and space variables 
\be
\tau=\epsilon^2 t, \qquad 
\color{black}
x= \alpha^{-1} \eta - ct,
\color{black}
\en 
where $c=\sqrt{\mu_0/\rho}$ is the speed of sound in the solid, and $\alpha$ is a suitable constant to be determined later. 
It is then possible to derive an $\mathcal{O}(\epsilon^3)$ asymptotic reduction  of \eqref{nw1}, as
\be \label{nw1bis}
F_\tau + \beta[(F^2+G^2)F]_x=0, \quad G_\tau + \beta[(F^2+G^2)C]_x=0,
\en
where $\beta$ is a  constant expressed in terms of $\rho$, $\mu_0$ and $\mu_1$.  
Vitolo and Saccomandi \cite{Vitolo}  discuss in  detail the mathematical and geometrical structures of the system \eqref{nw1} and of the corresponding asymptotic system  \eqref{nw1bis}. 
Note in particular that \eqref{nw1bis} is a Temple system, for which a general solution can be found  using a generalized hodograph transformation \cite{Tsarev}. 

Next we derive the asymptotic system of governing equations in the \emph{case of a pre-strain applied to a specific class of materials.}
Here we restrict attention to materials with a strain-energy density of \color{black} the form \cite{lopez2010new,destrade2017methodical}
\be \label{nw2}
\mathcal{W}=\frac{\mu_0}{2}	\left[I_1 - 3 + \frac{\kappa}{2}(I_1^2 - 3^2)\right],
\en
\color{black}
where $\mu_0>0$ is the infinitesimal shear modulus and $\kappa>0$ is a \color{black} measure of the departure from the new-Hookean model, in other words, a nonlinearity parameter with respect to wave propagation. \color{black} 
Then with $\hat F= \epsilon F$ and $\hat G=\epsilon G$ where $|\epsilon | \ll 1$, we differentiate the equations of motion \eqref{e7} with respect to $\eta$ to obtain the determining equations as 
\begin{align}
\label{nw3} 
&
 \rho \epsilon F_{tt} = [2\mathcal{W}_1\left(\vec{n \, \cdot} \overline{\vec{B}}\vec{a}+\epsilon F  \vec{n \, \cdot} \overline{\vec{B}}\vec{n}\right)]_{\eta \eta},
 \notag \\
& \rho \epsilon G_{tt}=[2\mathcal{W}_1\left(\vec{n \, \cdot} \overline{\vec{B}}\vec{b}+\epsilon G  \vec{n \, \cdot} \overline{\vec{B}}\vec{n}\right)]_{\eta \eta},
\end{align}
where
\begin{equation}
2 \mathcal{W}_1=\mu_0\left[1+\kappa \overline{I}_1 + 2\epsilon \kappa (F \vec{n \, \cdot} \overline{\vec{B}} \vec{a}+G \vec{n \, \cdot} \overline{\vec{B}} \vec{b})+\epsilon^2 \kappa \left(F^2+G^2 \right) \vec{n \, \cdot} \overline{\vec{B}} \vec{n}\right].
\end{equation}

We first consider that $\vec{n \, \cdot} \overline{\vec{B}} \vec{a}$ and $\vec{n \, \cdot} \overline{\vec{B}} \vec{b}$ are of order $\mathcal{O}(1)$ (for a detailed discussion about this point, see Pucci et al. \cite{Pucci}).
In this case the first nonlinearity arises at order $\epsilon^2$.
We introduce the scaled time and space variables
\begin{equation} \label{scaling}
\tau=\epsilon t, \qquad 
\color{black}
x = \alpha^{-1} \eta - ct,
\color{black}
\end{equation}
where $c= \sqrt{\mu_0 /\rho}$ is the speed of infinitesimal shear waves in the undeformed material and the non-dimensional constant $\alpha$ is to be determined soon. 
\color{black}
We then obtain the following equations,
\begin{equation} 
(\vec M - \alpha^2 \vec I)\begin{bmatrix} F\\G\end{bmatrix}_{xx} = 2\epsilon c^{-1}\begin{bmatrix} F\\G\end{bmatrix}_{\tau x} + 2 \epsilon \kappa( \vec{n \, \cdot} \overline{\vec{B}} \vec{n})  \left\{ \vec N \begin{bmatrix} F\\G\end{bmatrix}_{x}\right\}_x + \ldots,
\end{equation}
where $\vec M$ is a constant matrix,
\begin{equation} \label{M} 
\vec M = \begin{bmatrix}
(1+\kappa \overline{I}_1)  (\vec{n \, \cdot} \overline{\vec{B}}\vec{n}) + 2 \kappa  (\vec{n \, \cdot} \overline{\vec{B}} \vec{a})^2  &  2 \kappa (\vec{n \, \cdot} \overline{\vec{B}}\vec{a})( \vec{n \, \cdot} \overline{\vec{B}} \vec{b}) \\ 
 2 \kappa (\vec{n \, \cdot} \overline{\vec{B}}\vec{a})( \vec{n \, \cdot} \overline{\vec{B}} \vec{b})  & (1+\kappa \overline{I}_1)(  \vec{n \, \cdot} \overline{\vec{B}}\vec{n}) + 2 \kappa  (\vec{n \, \cdot} \overline{\vec{B}} \vec{b})^2
 \end{bmatrix},
\end{equation}
and $\vec N$ depends on $F$ and $G$,
\begin{equation} \label{N} 
\vec N = \begin{bmatrix}
3  (\vec{n \, \cdot} \overline{\vec{B}}\vec{a}) F +   (\vec{n \, \cdot} \overline{\vec{B}} \vec{b})G  &   (\vec{n \, \cdot} \overline{\vec{B}}\vec{b})F + ( \vec{n \, \cdot} \overline{\vec{B}} \vec{a}) G \\ 
 (\vec{n \, \cdot} \overline{\vec{B}}\vec{b}) F + ( \vec{n \, \cdot} \overline{\vec{B}} \vec{a}) G  & (\vec{n \, \cdot} \overline{\vec{B}}\vec{a}) F  + 3 (\vec{n \, \cdot} \overline{\vec{B}} \vec{b})G
 \end{bmatrix}.
\end{equation}

Expanding now $F$ and $G$ as
\begin{equation}
F = F_0(\tau,x) + \epsilon F_1(\tau,x) + \ldots, \qquad
G = G_0(\tau,x) + \epsilon G_1(\tau,x) + \ldots, 
\end{equation}
we then find the following equations at order $\epsilon^0$:
\begin{equation} \label{M0}
(\vec M - \alpha^2 \vec I)\begin{bmatrix} F_0\\G_0\end{bmatrix}_{xx} = \vec 0,
\end{equation}
and at order $\epsilon^1$:
\begin{equation} \label{N0}
(\vec M - \alpha^2 \vec I)\begin{bmatrix} F_1\\G_1\end{bmatrix}_{xx} = 
2 \left\{ c^{-1} \begin{bmatrix} F_0\\G_0 \end{bmatrix}_{\tau} + \kappa ( \vec{n \, \cdot} \overline{\vec{B}} \vec{n})   \vec N_0 \begin{bmatrix} F_0\\G_0\end{bmatrix}_{x}\right\}_x,
\end{equation}
where $\vec N_0$ is $\vec N$ in \eqref{N} when $F$, $G$ are replaced by $F_0$, $G_0$, respectively.

The system \eqref{M0} is an eigen-problem, with $\vec M$ symmetric. 
The eigenvalues are both positive:
\begin{align} \label{nw6} 
& \alpha_1^2 = (1+\kappa \overline{I}_1)(\vec{n \, \cdot} \overline{\vec{B}}\vec{n}),
\notag \\ &
  \alpha_2^2 = (1+\kappa \overline{I}_1)(\vec{n \, \cdot} \overline{\vec{B}}\vec{n}) +  2 \kappa [ (\vec{n \, \cdot} \overline{\vec{B}} \vec{a})^2+ (\vec{n \, \cdot} \overline{\vec{B}} \vec{b})^2],
\end{align}
with corresponding eigenvectors parallel to 
\begin{equation}
\vec v_1 = \begin{bmatrix} 
 \vec{n \, \cdot} \overline{\vec{B}}\vec{b}\\ -\vec{n \, \cdot} \overline{\vec{B}}\vec{a}
\end{bmatrix},
\qquad
 \vec v_2 = \begin{bmatrix} 
\vec{n \, \cdot} \overline{\vec{B}}\vec{a}\\ \vec{n \, \cdot} \overline{\vec{B}}\vec{b}
\end{bmatrix},
\end{equation}
respectively. 
It follows upon integration that $F_0$, $G_0$ are found as
\begin{equation}
\begin{bmatrix} F_0\\G_0\end{bmatrix} = \phi(\tau,x)
 \vec v_1 
\qquad 
\text{with} \quad 
x = \alpha_1^{-1}\eta - ct,
\end{equation}
or
\begin{equation}
\begin{bmatrix} F_0\\G_0\end{bmatrix} = \psi(\tau,x)
 \vec v_2
\qquad 
\text{with} \quad 
x =  \alpha_2^{-1}\eta - ct,
\end{equation}
depending on which value is chosen for $\alpha$ in the scaling \eqref{scaling}, where $\phi$, $\psi$ are arbitrary functions.

Taking now the dot product of \eqref{N0} with $\vec v_i$ ($i=1,2$), we arrive at
\begin{align} 
& \left[ (\vec v_1 \vec{\cdot v}_1) \phi_{\tau} + \kappa c ( \vec{n \, \cdot} \overline{\vec{B}} \vec{n})   (\vec v_1 \vec{\cdot }\vec N_0 \vec v_1) \phi_x \right]_x = 0, \notag \\
& \left[ (\vec v_2 \vec{\cdot v}_2) \psi_{\tau} + \kappa c ( \vec{n \, \cdot} \overline{\vec{B}} \vec{n})   (\vec v_2 \vec{\cdot }\vec N_0 \vec v_2) \psi_x \right]_x = 0, 
\end{align}
which we integrate with respect to $x$.
Simple calculations show that $\vec v_1 \vec{\cdot }\vec N_0 \vec v_1 = 0$, so that taking $\alpha = \alpha_1$ in the scaling \eqref{scaling} leads to a \emph{linearly degenerate wave} \cite{Hunter}.
On the other hand, $\vec v_2 \vec{\cdot }\vec N_0 \vec v_2 \ne 0$, and in that case we find a \emph{genuinely nonlinear wave} \cite{Hunter}, governed by the equation
\begin{equation}
  \psi_{\tau} + 3\kappa c ( \vec{n \, \cdot} \overline{\vec{B}} \vec{n})   [( \vec{n \, \cdot} \overline{\vec{B}} \vec{a})^2 + ( \vec{n \, \cdot} \overline{\vec{B}} \vec{b})^2]\psi   \psi_x  = 0.
  \end{equation}
Using the function $\bar \psi =  3\kappa c ( \vec{n \, \cdot} \overline{\vec{B}} \vec{n})   [( \vec{n \, \cdot} \overline{\vec{B}} \vec{a})^2 + ( \vec{n \, \cdot} \overline{\vec{B}} \vec{b})^2]\psi$, we see that this is the \emph{inviscid Burgers equation},
\begin{equation} \label{burgers}
\bar\psi_{\tau} + \left( \tfrac{1}{2} \bar  \psi^2\right)_x  = 0.
  \end{equation}

\color{black}

Finally, we consider the case where the direction of propagation of the wave is close to a principal direction, in the sense that $\vec{n \, \cdot} \overline{\vec{B}} \vec{a}$ and $\vec{n \, \cdot} \overline{\vec{B}} \vec{b}$ are of order $\mathcal{O}(\epsilon)$, \color{black} say $\vec{n \, \cdot} \overline{\vec{B}} \vec{a} = \epsilon a$ and $\vec{n \, \cdot} \overline{\vec{B}} \vec{b} = \epsilon b$, where $a$, $b$ are constants of order $\epsilon^0$. 

Then the scaling \eqref{scaling} does not lead to a consistent expansion. 
As pointed out by Pucci et al. \cite{Pucci} for the case of a single shear wave, the first nonlinearity arises at order $\epsilon^3$, as we must now use the  variables
\begin{equation}
\tau=\epsilon^2 t, \qquad 
x = \alpha^{-1} \eta - ct.
\color{black}
\end{equation}
Then we find that 
\begin{equation}
2 \mathcal{W}_1=\mu_0(1+\kappa \overline{I}_1 + \kappa \epsilon^2 \Lambda),
\end{equation}
where
\begin{equation} \label{Lambda}
\Lambda= 2 a F  + 2 b G  +  (\vec{n \, \cdot} \overline{\vec{B}} \vec{n})\left(F^2+G^2 \right).
\end{equation}

At order $\epsilon$ we obtain the decoupled system
\be \label{nw4bis} 
\alpha^2  F_{xx}=(1+\kappa \overline{I}_1)  (\vec{n \, \cdot} \overline{\vec{B}}\vec{n}) F_{xx}, \qquad  \alpha^2 G_{xx}=(1+\kappa \overline{I}_1) ( \vec{n \, \cdot} \overline{\vec{B}}\vec{n} ) G_{xx}
\en
which fixes  $\alpha$ as $\alpha = \sqrt{(1+\kappa \overline{I}_1)  (\vec{n \, \cdot} \overline{\vec{B}}\vec{n}) }$. 

The next terms are at order $\epsilon^3$, giving, upon integration, the system
\begin{align}  \label{nw8}
&  \alpha^2 F_\tau + \kappa c [(a + \vec{n \, \cdot} \overline{\vec{B}}\vec{n} \; F)\Lambda]_x=0, 
 \notag\\
&  \alpha^2 G_\tau + \kappa c [(b + \vec{n \, \cdot} \overline{\vec{B}}\vec{n} \; G)\Lambda]_x=0.
\end{align}
This coupled system contains both second- and third-order non-linearities.
It is consistent with the system \eqref{nw1bis} when $\vec{a}$ and $\vec{b}$ are aligned with the principal axes of $\overline{\vec{B}}$.
It can be written in the form
\begin{equation}
\alpha^2 \begin{bmatrix} F\\G\end{bmatrix}_{x} + \kappa c \; \vec A\begin{bmatrix} F\\G\end{bmatrix}_{\tau} =0,
\end{equation}
where
\begin{equation}
\vec A = \begin{bmatrix}
 (\vec{n \, \cdot} \overline{\vec{B}}\vec{n}) \Lambda +  2 [a + (\vec{n \, \cdot} \overline{\vec{B}} \vec{b})F]^2  &   2[a + (\vec{n \, \cdot} \overline{\vec{B}}\vec{n})F][b+ ( \vec{n \, \cdot} \overline{\vec{B}} \vec{n}) G] \\ 
   2[a + (\vec{n \, \cdot} \overline{\vec{B}}\vec{n})F][b+ ( \vec{n \, \cdot} \overline{\vec{B}} \vec{n}) G]
     &  (\vec{n \, \cdot} \overline{\vec{B}}\vec{n}) \Lambda +  2 [b + (\vec{n \, \cdot} \overline{\vec{B}} \vec{b})G]^2  
      \end{bmatrix},
\end{equation}
is symmetric and easy to diagonalise.
We find that the eigenvalues are
\begin{equation} \label{e-values}
\lambda_1 =  (\vec{n \, \cdot} \overline{\vec{B}}\vec{n}) \Lambda, \qquad
 \lambda_2= 2 (a^2  + b^2) + 3 (\vec{n \, \cdot} \overline{\vec{B}}\vec{n}) \Lambda,
\end{equation} 
with corresponding eigenvectors parallel to 
\begin{equation}
\begin{bmatrix}
 b+ ( \vec{n \, \cdot} \overline{\vec{B}} \vec{n}) G \\
-a - (\vec{n \, \cdot} \overline{\vec{B}}\vec{n})F
\end{bmatrix}, \qquad
\begin{bmatrix}
a+ ( \vec{n \, \cdot} \overline{\vec{B}} \vec{n}) F \\
b + (\vec{n \, \cdot} \overline{\vec{B}}\vec{n})G
\end{bmatrix},
\end{equation}
respectively. 
We find the Riemann invariants  as
\begin{equation} \label{Riemann}
R = \dfrac{a+ ( \vec{n \, \cdot} \overline{\vec{B}} \vec{n}) F}{b + (\vec{n \, \cdot} \overline{\vec{B}}\vec{n})G}, 
\qquad
S =  2 a F  + 2 b G  +  (\vec{n \, \cdot} \overline{\vec{B}} \vec{n})\left(F^2+G^2 \right).
\end{equation}
Therefore, bringing together \eqref{Lambda}, \eqref{e-values} and \eqref{Riemann}, we conclude that the eigenvalues  are written in terms of the Riemann invariants as
\begin{equation}
\lambda_1 =  (\vec{n \, \cdot} \overline{\vec{B}}\vec{n}) S, \qquad
 \lambda_2= 2 (a^2  +  b^2) + 3 (\vec{n \, \cdot} \overline{\vec{B}}\vec{n}) S,
\end{equation} 
and that the system is easily integrable in the hodograph plane.

\color{black}


\section{Concluding remarks}


Boulanger and Hayes were interested in the propagation of finite-amplitude shear waves in homogeneously deformed Mooney-Rivlin materials. 
They found that their motion is governed by linear differential equations and that, for a given propagation direction $\vec{n}$, only two directions of polarization are possible, those aligned with the principal axes of the elliptical section of the $\vec{x \cdot \bar B x}=1$ ellipsoid by the $\vec{n\cdot x} = 0$ plane. 
Then they turned their attention to the linear superposition of such two shear waves, in the process missing the generality that it is always possible to superpose two shear waves polarized in \emph{any} direction $ \vec{a} \perp \vec{n}$ (and $\vec b = \vec n \times \vec a$).

In fact, we found here that this result is valid for any isotropic strain-energy density function, not just the Mooney-Rivlin class. 
Moreover, we also complemented the results of the recent paper by Pucci et al. \cite{Pucci}.
The important findings in the full nonlinear setting are that the first nonlinearity that \emph{ matters} in pre-strained materials for shear waves is of second order, and that when the waves propagate and are polarized along, or close to, principal axes, (including the case $\overline{\vec{B}}=\vec{I}$), the first nonlinearity encountered is of third order. 

The advantage of using the asymptotic first-order equations \eqref{burgers} and \eqref{nw8} as compared to the full second-order equations is that they lend themselves to be solved using a plethora of analytical methods \cite{Dafermos}. 
We did not pursue this avenue here, as it was beyond the scope of this short note. 
We nonetheless  point out that adding a dissipative term to our constitutive equations, as in  \cite{Destrade2005}, would lead to a dissipative Burgers equation for \eqref{burgers} and a system of modified Burgers equations for \eqref{nw8}.
On the other hand, adding a dispersive term, as in \cite{DestradeSacco}, would turn  \eqref{burgers} into a KdV equation and the system \eqref{nw8}  into a set of coupled Gardner equations \cite{Gardner, GardnerB}. 
This is a remarkable fact, because Gardner's equation was introduced as a sort of mathematical `toy', and it is quite rare to encounter such an esoteric equation in the modelling of physics phenomena such as wave propagation.  
 

\begin{funding}

The research of MD was supported by a 111 Project for International Collaboration (Chinese Government, PR China) No. B21034 and by a grant from the Seagull Program (Zhejiang Province, PR China).
GS is supported by the Istituto Nazionale di Alta Matematica (INdAM), Fondi di Ricerca di Base UNIPG, the MIUR- PRIN project 2017KL4EF3 and Istituto Nazionale di Fisica Nucleare through its IS ‘Mathematical Methods in Non-Linear Physics’.

\end{funding}



\end{document}